\documentclass[11pt]{amsart}
\usepackage{mathrsfs}

\usepackage{amssymb}
\usepackage{amsfonts}

\usepackage{amscd}
\usepackage{bbm}

\usepackage[all]{xy}

\newtheorem{proposition}{Proposition}[section]

\newtheorem{definition}[proposition]{Definition}
\newtheorem{theorem}[proposition]{Theorem}

\newtheorem{conjecture}[proposition]{Conjecture}

\newtheorem{corollary}[proposition]{Corollary}

\newtheorem{remark}[proposition]{Remark}
\newtheorem{construction}[proposition]{Construction}
\newtheorem{lemma-definition}[proposition]{Lemma-Definition}

\newtheorem{mainaim}[proposition]{Main Goal}

\def\r{\color{red}}

\newcounter{tmp}

\def\coh{\operatorname{coh}}
\def\Qcoh{\operatorname{Qcoh}}

\def\lto{\longrightarrow}

\def\A{{\mathcal A}}

\def\D{{\mathcal D}}
\def\F{{\mathcal F}}
\def\E{{\mathcal E}}

\def\H{{\mathcal H}}

\def\N{{\mathcal N}}
\def\cO{{\mathcal O}}
\def\R{{\mathcal R}}
\def\L{{\mathcal L}}

\def\U{{\mathcal U}}

\def\M{{\mathcal M}}

\def\T{{\mathcal T}}

\def\ZZ{{\mathbb Z}}

\def\bD{{\mathbf D}}
\def\bR{{\mathbf R}}

\def\AA{{\mathbb A}}
\def\NN{{\mathbb N}}
\def\ZZ{{\mathbb Z}}

\def\PP{{\mathbb P}}

\def\Hom{\operatorname{Hom}}
\def\End{\operatorname{End}}
\def\Ext{\operatorname{Ext}}

\def\kk{{\mathbf k}}

\def\op{\circ}

\newcommand{\Ho}{{\H^0}}

\newcommand{\Ob}{\operatorname{Ob}}

\newcommand{\Ac}{\dA\!\mathit{c}\!\operatorname{--}\!}

\newcommand{\prf}{\mathcal{P}\!\mathit{erf}\!\operatorname{--}}

\def\dA{\mathscr A}
\def\dB{\mathscr B}
\def\dC{\mathscr C}

\def\dI{\mathscr I}

\def\Mod{{\mathscr M}\!\mathit{od}\!\operatorname{--}\!}

\def\mM{\mathsf M}

\def\mF{\mathsf F}

\def\mh{\mathsf h}

\def\mX{X}
\def\mY{Y}
\def\mZ{Z}

\def\dHom{\mathsf{Hom}}

\def\ptr{\operatorname{pre-tr}}

\def\wt{\widetilde}

\def\rd{\mathfrak{R}}

\newenvironment{dok}{\par\vspace{-5pt}%
\par\noindent\begingroup%
\leftskip=0em\hspace{0em}{\bf Proof.}}%
{\endgroup\hfill$\Box$}

{\endgroup\hfill$\Box$}

\def\a{A}

\def\wa{\widetilde{A}}

\def\db#1{ \bD^b({#1})}
\def\u{\mathbf u}
\def\bi{\mathbf i}
\def\r{\mathbf r}
\def\mod{\operatorname{mod}\!}

\def\wP{\widetilde{P}}
\def\wS{\widetilde{S}}
\def\wL{\widetilde{\L}}
\def\wX{\widetilde{X}}

\def\Hilb{\operatorname{Hilb}}

\title[]{Geometric realizations of quiver algebras}

\author[]{Dmitri Orlov}
\thanks{This work is supported by the Russian Science Foundation (RSF) under grant 14-50-00005}

\address{ Algebraic Geometry Department, Steklov Mathematical Institute of Russian Academy of Sciences,
Gubkin str. 8, Moscow 119991, RUSSIA}
\email{orlov@mi.ras.ru}

\date{}
\dedicatory{Dedicated to the blessed memory of my late adviser Andrei Nikolaevich Tyurin on the occasion of his 75th birthday}


\begin{document}
\begin{abstract}
In this paper we construct strong exceptional collections of vector bundles on smooth projective varieties 
that have a prescribed endomorphism algebra. We prove the construction problem always have a solution. 
We consider some applications to noncommutative projective planes and to the quiver connected with the 3-point Ising function.
\end{abstract}

\maketitle

\section*{Introduction}

The main purpose of this paper is to provide  constructions of strong exceptional collections on smooth projective varieties which consist of vector bundles and have a prescribed endomorphism algebra.
We are interested in triangulated categories $\T$ with a full exceptional collections $\sigma=(E_1,\dots, E_n).$
Recently, we showed that whenever such a $\T$ has an enhancement, i.e.
is equivalent to the homotopy category $\Ho(\dA)$ of some differential graded category $\dA,$  then it can be realized as an admissible subcategory of the bounded derived category of coherent sheaves on a smooth projective variety (see \cite[Th.5.8]{O_glue}).
Recall that a triangulated subcategory $\N\subset\mathbf{D}^b(\coh X),$ where $X$ is smooth and projective, is called admissible if it is full and the inclusion functor has right and  left adjoint projections.
Admissible subcategories have many good properties and provide a large selection of smooth and proper noncommutative schemes
that are called geometric noncommutative schemes \cite{O_glue}.

To provide some context we recall a result of \cite{O_glue} in more detail. Suppose the homotopy category
$\T=\Ho(\dA)$ of a small differential graded category $\dA$ has a full exceptional collection
$
\T=\langle E_1,\dots, E_n\rangle.
$
Then there exist a smooth projective scheme $X$ and an exceptional collection of line bundles
$\sigma=(\L_1,\dots, \L_n)$ on $X$ such that the full subcategory of $\db{\coh X},$
generated by $\sigma,$ is equivalent to $\T.$
In \cite{O_glue} we give an explicit construction of the variety $X$ as a tower of projective bundles.
It follows that $X$ itself
has a full exceptional collection.
Furthermore, we show that a full exceptional collection on $X$ can be chosen so that
it contains the collection $\sigma=(\L_1,\dots, \L_n)$ as a subcollection.
In this case we obtain a functor from the triangulated category $\T$ to the derived category
$\db{\coh X}$ that sends the exceptional objects $E_i$ to  shifts of the line bundles $\L_i[r_i]$ for some
integers $r_i.$ Of course, we can not expect in general that $E_i$ go to unshifted line bundles.

On the other hand, in the most important case when the exceptional collection
 $(E_1,\dots, E_n)$ is strong it is desirable to find a realizations of this collection
 in form of vector bundles (without shifts).
In this paper we will deal with strong exceptional collections and will discuss different constructions of geometric realizations of these collections  in term of vector bundles on smooth projective varieties. We show that
for a triangulated category $\T$ with a strong exceptional collection $\sigma=(E_1, \dots, E_n)$
it is always possible to find a smooth projective variety $X$ and a fully faithful functor
from $\T$ to the bounded derived category of coherent sheaves $\db{\coh X}$ such that
the exceptional objects $E_i$ go to vector bundles $\E_i$ on $X$ (see Theorem \ref{main} and Corollary \ref{exc_colllection_bundles}).
In this way, we obtain a strong exceptional collection $(\E_1, \dots, \E_n)$ of vector bundles
on $X$ with the same endomorphism algebra that the initial collection $\sigma$ has.

In the last section we consider some applications of our constructions to specific interesting exceptional
collections of three
objects. The first example is a quiver related to the  3-point Ising function, while the second example is the family of
quivers describing the noncommutative projective planes.

The author is very grateful to
Anton Fonarev, Alexander Kuznetsov, and Valery Lunts for very useful discussions and to Tony Pantev for a large number of valuable comments.

\section{Exceptional collections,  triangulated and differential graded categories}

\subsection{Exceptional collections}
We begin be recalling some definitions and facts concerning admissible
subcategories, semi-orthogonal decompositions, and exceptional collections (see \cite{BO}).
Let $\T$ be a $\kk$\!--linear triangulated category, where $\kk$ is
a base field. Let ${\N\subset\T}$ be a full triangulated
subcategory. Recall that the {\em right orthogonal} (resp. {\em left
orthogonal}) to ${\N}$ is the full
subcategory ${\N}^{\perp}\subset {\T}$ (resp. ${}^{\perp}{\N}$) consisting of all objects $X$
such that ${\Hom(Y, X)}=0$ (resp. ${\Hom(X, Y)}=0$)  for any $Y\in{\N}.$
It is clear that the
orthogonals are triangulated subcategories.

\begin{definition}\label{adm}
Let $j\colon\N\hookrightarrow\T$ be a full embedding of triangulated
categories. We say that ${\N}$ is {\em right admissible}
(resp. {\em left admissible}) if there is a right
(resp. left) adjoint functor $q\colon\T\to \N.$ The
subcategory $\N$ will be called {\em admissible} if it is both right
and left admissible.
\end{definition}

It is well-know that a subcategory $\N$ is right admissible if and only if for
each object $Z\in{\T}$ there is an exact triangle $Y\to Z\to X,$
with $Y\in{\N},\, X\in{\N}^{\perp}.$

Let $\N\subset\T$ be a full triangulated subcategory. If $\N$ is right (resp. left) admissible, then
the quotient category $\T/\N$ is equivalent to $\N^{\perp}$
(resp. ${}^{\perp}\N$).  Conversely, if the quotient functor
$\T\lto\T/\N$ has a left (resp. right) adjoint,
then $\T/\N$ is equivalent to $\N^{\perp}$
(resp. ${}^{\perp}\N$).

\begin{definition}\label{sd}
A {\em semi-orthogonal decomposition} of a triangulated category
$\T$ is a sequence of full triangulated subcategories ${\N}_1,
\dots, {\N}_n$ in ${\T}$ such that there is an increasing filtration
$0=\T_0\subset\T_1\subset\cdots\subset\T_n=\T$ by left admissible
subcategories for which the left orthogonals ${}^{\perp}\T_{p-1}$ in
$\T_{p}$ coincide with $\N_p.$ In particular,
$\N_p\cong\T_p/\T_{p-1}.$ We write $ {\T}=\left\langle{\N}_1, \dots,
{\N}_n\right\rangle. $
\end{definition}

In some cases one can hope that $\T$ has a semi-orthogonal decomposition
${\T}=\left\langle{\N}_1, \dots, {\N}_n\right\rangle$ in which each
$\N_p$ is as simple as possible, i.e. each of them is equivalent to the bounded
derived category of finite-dimensional vector spaces.

\begin{definition}\label{exc}
An object $E$ of a $\kk$\!--linear triangulated category ${\T}$ is
called {\em exceptional} if  ${\Hom}(E, E[l])=0$ whenever $l\ne 0,$
and ${\Hom}(E, E)=\kk.$ An {\em exceptional collection} in ${\T}$ is
a sequence of exceptional objects $\sigma=(E_1,\dots, E_n)$
satisfying the semi-orthogonality condition ${\Hom}(E_i, E_j[l])=0$
for all $l$ whenever $i>j.$
\end{definition}

If a triangulated category $\T$ has an exceptional collection
$\sigma=(E_1,\dots, E_n)$ that generates the whole of $\T,$ then
this collection is called {\em full}.  In this case $\T$ has a
semi-orthogonal decomposition with $\N_p=\langle E_p\rangle.$ Since
$E_{p}$ is exceptional, each of these categories is equivalent to
the bounded derived category of finite dimensional $\kk$\!-vector
spaces.  In this case we write $ \T=\langle E_1,\dots, E_n \rangle.$

\begin{definition}\label{strong}
An exceptional collection $\sigma=(E_1,\dots, E_n)$ is called {\em
strong} if, in addition, ${\Hom}(E_i, E_j[l])=0$ for all  $i$ and
$j$ when $l\ne 0.$
\end{definition}

Let $\T$ be a triangulated category with a full strong exceptional
collection $\sigma=(E_1,\dots, E_n).$ The algebra $\a=\End
(\mathop{\bigoplus}\limits_{i=1}^n E_i)$
is called the {\em algebra of endomorphisms of the exceptional collection $\sigma.$}

\subsection{Differential graded categories and enhancements}
Here we only introduce notations and recall some facts on differential graded (DG) categories.
Our main references for DG categories are \cite{Ke,Dr} (see also \cite{LO, O_glue}).
A {\em differential graded or DG category} is a $\kk$\!--linear category $\dA$ whose morphism spaces
$\Hom (\mX, \mY)$
are complexes of $\kk$\!-vector spaces, so that for any $\mX, \mY, \mZ\in
\Ob\dC$ the composition $\Hom (\mY, \mZ)\otimes \Hom (\mX, \mY)\to
\Hom (\mX, \mZ)$ is a morphism of DG $\kk$\!--modules.
The identity morphism $1_\mX\in \Hom (\mX, \mX)$ is required to be closed of
degree zero.

For a DG category $\dA$ we denote by $\Ho(\dA)$
its homotopy category. The homotopy category $\Ho(\dA)$ has the same objects as the DG category $\dA$ and its
morphisms are defined by taking the $0$\!-th cohomology
$H^0\Hom_{\dA} (\mX, \mY)$
of the complex $\Hom_{\dA} (\mX, \mY).$

As usual, a DG functor
$\mF:\dA\to\dB$ is given by a map $\mF:\Ob(\dA)\to\Ob(\dB)$ and
by morphisms of DG $\kk$\!--modules
$$
\mF_{\mX, \mY}: \dHom_{\dA}(\mX, \mY) \lto \dHom_{\dB}(\mF \mX,\mF \mY),\quad \mX, \mY\in\Ob(\dA)
$$
compatible with the composition and units.
A DG functor $\mF: \dA\to\dB$ is called a {\em quasi-equivalence} if
$\mF_{\mX, \mY}$ is a quasi-isomorphism for all pairs of objects $\mX, \mY$ of $\dA$
and the induced functor $H^0(\mF): \Ho(\dA)\to \Ho(\dB)$ is an
equivalence. Two DG categories $\dA$ and $\dB$ are {\em quasi-equivalent} if there is a DG
category $\dC$ and quasi-equivalences
$\dA\stackrel{\sim}{\leftarrow} \dC \stackrel{\sim}{\rightarrow}\dB.$

Given a small DG category $\dA$ we define a right DG $\dA$\!--module as a DG functor
$\mM: \dA^{op}\to \Mod \kk,$ where $\Mod \kk$ is the DG category of DG $\kk$\!--modules. We denote by $\Mod \dA$ the DG
category of right DG $\dA$\!--modules.
Each object $\mY$ of $\dA$ produces a right module
$
\dHom_{\dA}(-, \mY)
$
which is called a free DG module represented by $\mY.$ We obtain the Yoneda DG functor
$\mh^\bullet :\dA \to
\Mod\dA$ that is fully faithful.
The {\it derived
category} $\D(\dA)$ is defined as the Verdier quotient
\[
\D(\dA):=\Ho(\Mod\dA)/\Ho (\Ac\dA),
\]
where $\Ac\dA$ is the full
DG subcategory of $\Mod\dA$ consisting of all acyclic DG modules, i.e. DG modules $\mM$
for which the complexes of $\kk$\!-modules $\mM(\mX)$ are acyclic for all $X\in\dA.$

\begin{definition}
The triangulated category of {\em perfect DG modules} $\prf \dA$
is the smallest triangulated subcategory of $\D(\dA)$ that contains all free DG modules and is closed under direct summands.
\end{definition}

The triangulated categories $\D(\dA)$ and $\prf\dA$ are invariant under quasi-equivalences of 
$\dA.$

For any DG category $\dA$ there exist a DG category $\dA^{\ptr}$ that is called the pretriangulated hull
of $\dA$
and a canonical fully faithful DG
functor $\dA\hookrightarrow\dA^{\ptr}.$
The idea of the definition of $\dA^{\ptr}$ is to add to $\dA$
all shifts, all cones, cones of morphisms between cones and etc.
A DG category $\dA$ is called pretriangulated  if the canonical DG functor $\dA\to\dA^{\ptr}$
is a quasi-equivalence.
 It is equivalent to require that the homotopy category $\Ho(\dA)$ is triangulated
as a subcategory of $\Ho(\Mod\dA).$
The DG category $\dA^{\ptr}$ is always  pretriangulated,
so $\Ho(\dA^{\ptr})$ is a triangulated category.

\begin{definition} Let $\T$ be a triangulated category. An {\em
enhancement} of $\T$ is a pair $(\dA , \varepsilon),$ where $\dA$ is a
pretriangulated DG category and $\varepsilon:\Ho(\dA)\stackrel{\sim}{\to} \T$ is an exact equivalence.
\end{definition}

For any quasi-compact and separated  scheme $X$ over an arbitrary field $\kk$
the derived category $\bD(\Qcoh X)$ has an enhancement
that is coming from h-injective complexes (see, e.g. \cite{KSh}), i.e.  $\Ho(\dI(X))\cong \bD(\Qcoh X),$ where $\dI(X)$ the DG category
of h-injective complexes. As a consequence, we obtain an enhancements for any full triangulated subcategory of $\bD(\Qcoh X),$
for example for the triangulated  category of perfect complexes $\prf X$ and for the bounded derived category of coherent sheaves  $\db{\coh X}$ in noetherian case.

There are different notions of generators in triangulated categories.
We recall the most useful notion of generating
a triangulated category that is the notion of a classical generator.

\begin{definition}
An object $E\in \T$
is called a {\em  classical generator} if the category $\T$ coincides with the smallest triangulated subcategory
that contains $E$ and is closed under taking direct summands.
\end{definition}
 Note that the category of perfect complexes $\prf X$ admits a classical generator for any quasi-compact and quasi-separated scheme $X$ (see \cite{Ne, BV}). If $X$ is quasi-projective of dimension $d$ then the object $\mathop{\bigoplus}\limits_{p=0}^{d}\L^{\otimes p},$ where $\L$ is very ample line bundle, is a classical generator (see \cite{O_gen}).

The following theorem shows that the notion of classical generator is very useful for triangulated categories admitting enhancements.

\begin{theorem}[\cite{Ke, Ke2}] \label{keller}
Let $\T$ be an idempotent complete triangulated category that admits an enhancement $\dA$ and let $E\in \T$
be a classical generator. Then the category $\T$ is equivalent to the triangulated category $\prf \A,$
where $\A=\dHom_{\dA}(E, E)$ is the DG algebra of endomorphisms of the object $E$ in the DG category $\dA.$
\end{theorem}

This theorem implies the following corollary that we will use in the sequel.

\begin{corollary}\label{exceptional_generator}
Let $\T$ be a triangulated category that admits an enhancement.
Assume that $\T$ has a full strong exceptional collection $\sigma=(E_1,\dots, E_n).$
Then the category $\T$ is equivalent to the derived category $\db{\mod-\a},$ where $\a=\End(\mathop{\bigoplus}\limits_{i=1}^n E_i)$ is the algebra of endomorphisms of the collection $\sigma.$
\end{corollary}
\begin{dok}
Since $\sigma$ is full, the object $E=\mathop{\bigoplus}\limits_{i=1}^n E_i$ is a classical generator. As $\sigma$ is strong,
the DG algebra of endomorphisms of the object $E$ has only 0-th cohomology. Hence, this DG algebra is quasi-isomorphic to the usual
endomorphism algebra of the collection $\sigma.$ As a triangulated category with a full exceptional collection
 the category $\T$ is idempotent complete. Now corollary follows from the Theorem \ref{keller} and the fact
that for any algebra of finite global dimension the derived category $\db{\mod-\a}$ is equivalent to the category of perfect complexes over $\a.$
\end{dok}

In the paper \cite{O_glue} we showed that any triangulated category with a full exceptional collection has a geometric realization
as long as it has an enhancement. More precisely,
we proved the following.
\begin{theorem}\label{exc_col}\cite[Th.5.8]{O_glue}
Let $\dA$ be a small DG category over $\kk$ such that the homotopy category
$\T=\Ho(\dA)$ has a full exceptional collection
$
\T=\langle E_1,\dots, E_n\rangle.
$
Then there are a smooth projective scheme $X$ and an exceptional collection of line bundles
$\sigma=(\L_1,\dots, \L_n)$ on $X$ such that the subcategory of $\db{\coh X},$
generated by $\sigma,$ is equivalent to $\T.$
Moreover, $X$ is a sequence of projective bundles and has a full exceptional collection.
\end{theorem}

The scheme $X$ has a full exceptional collection as a tower of projective bundles
(see \cite{Blow}). Furthermore, it follows from the construction that a full exceptional collection on $X$ can be chosen in a way that
it contains the collection $\sigma=(\L_1,\dots, \L_n)$ as a subcollection.

In the proof of this theorem we constructed  a functor from the triangulated category $\T$ to the derived category
$\db{\coh X}$ that sends the exceptional objects $E_i$ to shifts of the line bundles $\L_i[r_i]$ for some
integers $r_i.$ Of course, we can not expect in general that $E_i$ go to line bundles without shifts.
On the other hand, in the case of  strong exceptional collections it is natural to seek realizations as 
 collections of vector bundles (without shifts) on  smooth projective varieties. It can be shown that in general
 we can not realize a strong exceptional collection as a collection of unshifted line bundles (see  Remark \ref{remaark_line}),
but trying to present it in terms of vector bundles seem quite reasonable.

In this paper we deal with strong exceptional collections and  discuss different constructions of geometric realizations 
as strong exceptional collections  of vector bundles on smooth projective varieties. We prove that such realizations  always exist.

\section{Geometric realizations}

\subsection{Quiver algebras}
 A quiver is a finite directed graph, possibly with multiple arrows
and loops. More precisely, a quiver $Q$ consists of a data
$(Q_0, Q_1, s, t),$
where $Q_0, Q_1$ are finite sets of vertices and arrows respectively, while
$s, t : Q_1\to  Q_0$
are maps attaching to each arrow its source and target.

The path $\kk$\!-algebra of the quiver $Q$ is the algebra $\kk Q$ determined
by the generators $e_q$ for $q\in Q_0$ and $a$ for $a\in Q_1$
with the following  relations
\[
e^2_q = e_q,\quad e_r e_q = 0, \quad\text{ when}\quad r\ne q,
\quad\text{and}
\quad
e_{t(a)} a = a e_{s(a)} = a.
\]
In particular, the elements $e_q$ are orthogonal idempotents of
the path algebra $\kk Q.$
It follows from the relations above that
$e_q a = 0$ unless $q=t(a)$ and $a e_q = 0$ unless $q = s(a).$

As a $\kk$\!-vector space
the path algebra $\kk Q$ has a basis consisting of the set of all paths in $Q,$
where a path $\overline{p}$ is a possibly empty sequence $a_{m} a_{m-1}\cdots a_1$ of compatible arrows,
i.e. $s(a_{i+1})=t(a_i)$ for all $i=1,\dots, m-1.$
For an empty path we have to choose a vertex from $Q_0.$
The composition of two paths $\overline{p}_1$ and $\overline{p}_2$ in $Q$
is defined naturally as $\overline{p}_2 \overline{p}_1$ if they are compatible
and as $0$ if they are not compatible.
This is a more natural definition of  the
product in paths algebra $\kk Q.$

To obtain a more general class of algebras, it is useful to introduce in consideration  quivers with
relations.
A relation on a quiver $Q$ is a subspace of $\kk Q$ spanned by linear
combinations of paths having a common source and a common target, and of length at
least 2.
A quiver with relations is a pair $(Q, I),$ where $Q$ is a quiver and $I$ is a two-sided ideal
of the path algebra $\kk Q$ generated by relations.
The quotient algebra $\kk Q/I$ will be called a quiver algebra of the quiver with relations
$(Q, I).$
It can be shown that every module
category $\mod-\a,$  where $\a$ is a finite dimensional algebra over $\kk,$ is equivalent
to $\mod-\kk Q/I$ for some quiver with relations $(Q, I)$(see e.g. \cite{Gab}).

A quiver algebra $\a=\kk Q/I$ viewed  as right module over itself can be decomposed as a direct sum
of projective modules $P_q=e_q\a$ for $q\in Q_0,$ i.e. $A=\mathop{\oplus}\limits_{q\in Q_0} P_q.$
The projective modules $P_q\subset \a$ consist
of all paths $\overline{p}$ with fixed target (or tail) $t(\overline{p})=q.$

In this paper we consider quiver algebras for special type quivers with relations
that are directly related to
exceptional collections.

\begin{definition}
We say that $\a$ is a {\em quiver algebra on $n$ ordered vertices} if it is a quiver algebra of
a quiver with relations
$(Q, I)$ for which $Q_0=\{1,\dots , n\}$ is the ordered set of $n$ elements  and  $s(a)< t(a)$
for any arrow $a\in Q_1.$
\end{definition}

It is evident that the algebra of endomorphisms of any (strong) exceptional collection
$\sigma=(E_1, \dots E_n)$  is a quiver algebra
on $n$ ordered vertices.
On the other hand, a quiver algebra $\a$ on $n$ ordered vertices has finite global dimension and,
moreover, its derived category $\db{\mod-\a}$ has a strong full exceptional collection
consisting of the projective modules $P_i$ for $i=1,\dots, n.$ The algebra $\a$ is exactly
the algebra of endomorphisms of this full strong exceptional collection.

\subsection{Strong exceptional collections and geometric realizations}

Let $\a$ be a quiver algebra on $n$ ordered vertices.
Denote by $P_1,\dots, P_{n}$ the respective right projective modules.
As it was mentioned above, the collection of projective modules $(P_1,\dots, P_{n})$ is a full strong exceptional collection in the category $\db{\mod-\a}.$

\begin{mainaim} Let $\a$ be a quiver algebra on $n$ ordered vertices. Our main goal is to construct a smooth projective variety
$X$ and a strong (not full) exceptional collection $\sigma=(\E_1,\dots, \E_n)$ of vector bundles on $X$ such that the algebra
of endomorphisms $\End(\mathop{\bigoplus}\limits_{i=1}^{n} \E_i)$ coincides with the algebra  $\a.$
\end{mainaim}

Assume that on a smooth projective variety
$X$ there is a strong exceptional collection $\sigma=(\E_1,\dots, \E_n),$ for which the algebra
of endomorphisms $\End(\mathop{\bigoplus}\limits_{i=1}^{n} \E_i)$ coincides with the algebra  $\a.$
In this case by Corollary \ref{exceptional_generator}
the admissible subcategory of $\db{\coh X}$ generated by the collection $\sigma$
is equivalent to the derived category
$\db{\mod-\a}$ of the endomorphism algebra $\a.$
The inclusion functor sends
the projective modules $P_i$ to the objects $\E_i.$
The right adjoint projection functor is
\[
\bR\Hom^{\cdot}(\E, -): \db{\coh X}\lto \db{\mod-\a},
\quad
\text{where}
\quad
\E=\mathop{\bigoplus}\limits_{i=1}^{n} \E_i.
\]
All of these statements are true for any objects $\E_i$ not only for vector bundles.
And it is proved in \cite{O_glue} that such realizations exist for any exceptional collection.
Our main aim in this paper is to give constructions of  strong exceptional collections
in terms of vector bundles.
The main technique that we will use is an induction on the number of vertices and a passage from a quiver algebra $\a$
to its ordinary extension by an $\a$\!-module.

Let $\a$ be a quiver algebra on $n$ ordered vertices. Denote by $S_i, i=1,\dots, n,$ its respective simple modules.
The set of simple modules $(S_n, \dots, S_1)$ also forms a full (not strong) exceptional collection in $\db{\mod-\a}$  called the dual to the collection of projective modules.
Thus we have the following semi-orthogonal decompositions
\[
\db{\mod-\a}=\left\langle P_{1}, \dots, P_{n} \right\rangle=\left\langle S_{n}, \dots, S_1 \right\rangle.
\]

Note that for any  finite right $\a$\!-module $M$ there is
a filtration $0=M_0\subset M_1\subset \cdots\subset M_n=M$ such that each successive quotient
$M_p/M_{p-1}$ is isomorphic
to a (finite) direct sum of copies of the corresponding simple module $S_p.$

Let us consider  a quiver algebra $\wa$ on $n+1$ ordered vertices.
Denote by $\wP_1,\dots, \wP_{n+1}$ its right projective modules and
denote by $\wS_i$ the respective simple modules.

Consider the first $n$ vertices and denote by $\a$ the quiver algebra of the subquiver
corresponding to  these $n$ vertices.
We have a full embedding $\mod-\a\to\mod-\wa$ of abelian categories. It is an exact functor that sends the simple $\a$\!-modules $S_i$ to the simple $\wa$\!-modules $\wS_i$ and the projective $\a$\!-modules
$P_i$ to the projective $\wa$\!-modules $\wP_i$ for all $i=1,\dots,n.$
Furthermore,
the induced derived functor
$\bi: \db{\mod-\a}\to\db{\mod-\wa}$ is fully faithful.

Thus, there are semi-orthogonal decompositions
\[
\db{\mod-\wa}=\left\langle \wS_{n+1},\; \db{\mod-\a} \right\rangle=\left\langle \db{\mod-\a},\; \wP_{n+1} \right\rangle.
\]
Let us consider the natural surjection $\wP_{n+1}\twoheadrightarrow \wS_{n+1}$ and denote by $M$ the kernel
of this map. In the exact sequence of $\wa$\!-modules
\[
0\lto M\lto \wP_{n+1}\lto \wS_{n+1}\lto 0
\]
the module $M$ belongs to the subcategory $\bi( \db{\mod-\a})$ and, hence,
it can be considered as $\a$\!-module.
In particular, for all $i=1,\dots, n$ we have isomorphisms
\[
\Hom_{\wa}(\wP_i,\; \wP_{n+1})\cong \Hom_{\wa}(\wP_i,\; M)\cong
\Hom_{\a}(P_i,\; M).
\]
This means that the algebra $\wa$ can be constructed as low-triangular algebra of the following form
\begin{equation}\label{simplext}
\wa=\begin{pmatrix}
A & 0\\
M & \kk
\end{pmatrix},
\end{equation}
where $\a$ is a quiver algebra on $n$ vertices and $M$ is a right $\a$\!-module. The algebra $\wa$ is uniquely determined by
the algebra $\a$ and the right $\a$\!-module $M.$

\begin{definition}
The algebra $\wa$ defined by rule (\ref{simplext}) will be called an  {\em ordinary extension}
of the algebra $\a$ via the $\a$\!-module $M.$
\end{definition}

We start with a realization of the algebra $\a$ in terms of vector bundles and will try to extend this to a realization on the algebra $\wa.$
First, we record a simple fact that will be useful for applications.

\begin{proposition}\label{simple} Let $\a$ be a quiver algebra on $n$ ordered vertices and let $\wa$ be its ordinary extension via an $\a$\!-module $M.$
Suppose there are a smooth projective variety $X$ and a fully faithful functor
$\r: \db{\mod-\a}\to \db{\coh X}$ that sends the projective
modules $P_i$ to vector bundles $\E_i.$
Assume that there is a vector bundle $\M$ on $X$ such that
$\bR\Hom^{\cdot}(\mathop{\oplus}\limits_{i=1}^{n} \E_i,\; \M)\cong M$
in $\db{\mod-\a}.$
Then there are a smooth projective variety $\wX$
and a fully faithful functor
$\wt{\r}: \db{\mod-\wa}\to \db{\coh \wX}$ that sends the projective
modules $\wP_i$ to some vector bundles $\wt{\E}_i$ on $\wX.$
\end{proposition}
\begin{dok}
The proof is direct. Let us put $\wX=\PP(\M^{\vee})$ with the natural projection $\pi: \wX\to X.$
There is a canonical exact sequence on $\wX$
\begin{equation}\label{grothendieck}
0\lto \Omega_{\wX/X}(1)\lto \pi^*\M \lto \cO_{\wX}(1)\lto 0,
\end{equation}
where $\Omega_{\wX/X}$ is the relative cotangent bundle and $\cO_{\wX}(-1)$ is the tautological line bundle on $\wX.$

Let us take $\wt{\E}_i=\pi^*(\E_i)$ for all $i=1, \dots, n$ and put $\wt{\E}_{n+1}=\cO_{\wX}(1).$
It is easy to see that the collection $\wt{\sigma}=(\wt{\E}_1, \dots, \wt{\E}_n, \wt{\E}_{n+1})$
is exceptional.
There is a sequence of  isomorphisms
\[
\Ext^j_{\wX}(\wt{\E}_i,\; \wt{\E}_{n+1})=
\Ext^j_{\wX}(\pi^*\E_i,\; \cO_{\wX}(1))\cong
\Ext^j_{X}(\E_i,\; \bR\pi_*\cO_{\wX}(1))\cong
\Ext^j_{X}(\E_i,\; \M)\cong
\Ext^j_{A}(P_i,\; M)
\]
for all $i=1,\dots, n.$ It implies that the exceptional collection $\wt{\sigma}$ is strong and the algebra of endomorphisms of this collection
is isomorphic to the extended algebra $\wa.$

There is another construction when one takes as $\wX$ the projective bundle $\PP(\M).$ In this case
the dual exact sequence should be considered
\begin{equation}\label{grothendieck2}
0\lto \cO_{\wX}(-1)\lto \pi^*\M \lto \T_{\wX/X}(-1)\lto 0.
\end{equation}
As above we take $\wt{\E}_i=\pi^*(\E_i)$ for all $i=1, \dots, n,$ but put $\wt{\E}_{n+1}=\T_{\wX/X}(-1).$
The similar calculations give us that the collection $\wt{\sigma}=(\wt{\E}_1, \dots, \wt{\E}_n, \wt{\E}_{n+1})$
is strong exceptional and its algebra of endomorphisms
is also isomorphic to the ordinary extended algebra $\wa.$

In both cases we obtain a fully faithful functor $\wt{\r}:
\db{\mod-\wa}\to \db{\coh \wX}$ that sends the projective
modules $\wP_i$ to the vector bundles $\wt{\E}_i$ on $\wX$  for all $i=1,\dots, n+1.$
\end{dok}

These constructions can be useful in particular cases when one would like to represent a strong exceptional collection
as a collection of vector bundles on a smooth projective variety.
 However, they can not help with the general proof, because we can not guarantee that
 a given module will be represented by a vector bundle in the next step of induction.

\subsection{The main theorem}

Let $\a$ be a quiver algebra on $n$ ordered vertices.
Suppose there is an exact functor  $u: \mod-A\to \coh(X)$ between abelian categories.
Denote by $\u$ the derived functor from $\db{\mod-\a}$ to $\db{\coh X}.$ Since $u$ is exact the derived functor
 is defined and coincides with $u$ on the abelian category $\mod-A,$ i.e. $\u(M)\cong u(M)$ for any $\a$\!-module $M.$

We say that the exact functor $u: \mod-A\to \coh(X)$ satisfies property (V) if the following conditions hold
\begin{equation*}\label{V}
\parbox{.93\textwidth}{\begin{enumerate}
\item[1)] the induced derived functor $\u:\db{\mod-\a}\to \db{\coh X}$ is fully faithful;
\item[2)] simple modules $S_i$  go to line bundles $\L_i$ on $X$ under $u;$
\item[3)] there is a line bundle $\L$ on $X$ such that $\L\in {}^{\perp}\u(\db{\mod-\a}),$ the line bundles
$\L\otimes \L_i^{-1}$ are generated by global sections  and
$H^j(X, \L\otimes \L_i^{-1})=0$ when $j\ge 1$ for all $i=1,\dots,n.$
\end{enumerate}
}\leqno {\rm (V)}
\end{equation*}

Since any module $M$ has a composition filtration with successive quotient being simple modules,
the condition (2) implies that any $\a$\!-module $M$ goes to a vector bundle under the functor $u.$
Moreover, by the same reasoning the vector bundle $u(M)$ has a filtration with successive quotients isomorphic to the line bundles
$\L_i.$ Now it is not difficult to check that  condition (3) of (V) implies the following condition
\begin{equation*}
\parbox{.88\textwidth}{\begin{enumerate}
\item[3')] there is a line bundle $\L$ on $X$ such that $\L\in {}^{\perp}\u(\db{\mod-\a})$ and the vector bundles
$\L\otimes u(M)^{\vee}$ are generated by global sections for all $M\in\mod-\a.$
\end{enumerate}}
\end{equation*}

Consider the exact functor $u: \mod-A\to \coh(X)$ and denote by $\E_i$ the vector bundles $u(P_i)$ for all $i=1,\dots, n.$
Since the derived functor $\u$ is fully faithful, the sequence $\sigma=(\E_1,\dots, \E_n)$ of vector bundles on $X$
forms a strong exceptional collection.

\begin{proposition}\label{step} Let $\a$ be a quiver algebra on $n$ ordered vertices and let $\wa$ be its ordinary extension via an $\a$\!-module $M.$
Suppose there exist a smooth projective scheme $X$ and
 an exact functor
$u:\mod-\a\to \coh(X)$ that satisfies property (V).
Then there are a vector bundle $\F$ on $X$ and an exact functor $\wt{u}: \mod-\wa\to \coh(\PP(\F))$ that also satisfies
property (V). Moreover, the restriction of the functor $\wt{u}$ on $\mod-\a$ is isomorphic to $\pi^*\circ u,$
where $\pi:\PP(\F) \to X$ is the natural projection.
\end{proposition}

\begin{dok}
Let $u:\mod-\a\to \coh(X)$ be an exact fully faithful functor that satisfies property (V).
Denote by $\M$ the vector bundle $u(M).$ By assumption (3) and its consequence (3') there is a surjection
$(\L^{-1})^{\oplus m}\twoheadrightarrow\M^{\vee}$ for some $m\in \NN.$
Denote by $\F$ the vector bundle on $X$ that is dual to the kernel of this surjection.
Thus, we have an exact sequence
\begin{equation}\label{ext}
0\lto \F^{\vee}\lto (\L^{-1})^{\oplus m}\lto\M^{\vee}\lto 0.
\end{equation}

Taking sufficiently large $m$ we can assume that the rank of $\F$ is greater  than $2.$
Let us consider the projective bundle $\pi: \PP(\F)\to X$ and denote it as $\wX.$
There are natural exact sequences on $\wX$ of the following form
\[
0\lto \Omega_{\wX/X}(1)\lto \pi^*\F^{\vee}\lto \cO_{\wX}(1)\lto 0,\quad\text{and}\quad
0\lto \cO_{\wX}(-1)\lto \pi^*\F\lto \T_{\wX/X}(-1)\lto 0,
\]
where $\cO_{\wX}(-1)$ is the tautological line bundle and $\T_{\wX/X}, \Omega_{\wX/X}$ are relative tangent and relative cotangent bundles, respectively.

Denote by  $\wt{\L}_i$ the pull back line bundles $\pi^*\L_i$ for $i=1,\dots, n$ and put  $\wL_{n+1}=\cO_{\wX}(-1).$
Under the sequence of isomorphisms
\[
\Ext^1_X(\M^{\vee},\; \F^{\vee})\cong \Ext^1_{\wX}(\pi^*\M^{\vee},\; \cO_{\wX}(1))\cong\Ext^1_{\wX}(\cO_{\wX}(-1),\; \pi^*\M)
\]
the element $e\in\Ext^1_X(\M^{\vee},\; \F^{\vee}),$ which defines the short exact sequence (\ref{ext}),
gives some element $e'\in\Ext^1_{\wX}(\cO(-1),\; \pi^* \M).$ The element $e'$ defines the following extension
\begin{equation}\label{projective}
0\lto\pi^*\M\lto \wt{\E}_{n+1}\lto\cO(-1)\lto 0.
\end{equation}
that can be considered as a definition of the vector bundle $\wt{\E}_{n+1}.$

It follows from the definition of the exact sequence (\ref{projective}) that  the dual sequence goes to the the exact sequence
(\ref{ext}) under the direct image functor $\bR\pi_*.$
In particular, there is an isomorphism
$\bR\pi_* \wt{\E}_{n+1}^{\vee}\cong (\L^{-1})^{\oplus m}.$
This fact implies the following
vanishing of Hom's spaces
\begin{multline}\label{vanish}
\Hom^{i}_{\wX}(\wt{\E}_{n+1},\; \pi^*\N)\cong
H^i(\wX,\; \pi^*\N\otimes\wt{\E}_{n+1}^{\vee})\cong
H^{i}(X,\; \N\otimes\bR\pi_* \wt{\E}_{n+1}^{\vee})
\cong\\
 \cong H^{i}(X,\; \N\otimes (\L^{-1})^{\oplus m})
\cong \Hom^{i}_X(\L^{\oplus m},\; \N)=0
\end{multline}
for any object $\N$ from the image of the functor $\u$ by (3) of (V).
Hence, the vector bundle $\wt{\E}_{n+1}$ belongs to the left orthogonal ${}^{\perp}\pi^*\u (\db{\mod-A}).$

Furthermore, since $\Hom^{\cdot}_{\wX}(\wt{\E}_{n+1},\; \pi^*\M)=0$ and
$\Hom^{\cdot}_{\wX}(\pi^*\M,\; \cO_{\wX}(-1))=0,$ the short exact sequence (\ref{projective}) induces the following isomorphisms
\begin{equation}\label{except}
\Hom^{\cdot}_{\wX}(\wt{\E}_{n+1},\; \wt{\E}_{n+1})\cong
\Hom^{\cdot}_{\wX}(\wt{\E}_{n+1},\; \cO_{\wX}(-1))\cong
\Hom^{\cdot}_{\wX}(\cO_{\wX}(-1),\; \cO_{\wX}(-1)).
\end{equation}
This implies that $\wt{\E}_{n+1}$ is exceptional.
Denote by $\wt{\E}_{i}$ the vector bundles $\pi^*\E_{i}$ for all $i=1,\dots, n.$
The vanishing properties (\ref{vanish}) and the isomorphisms (\ref{except}) imply
that the ordered set  $(\wt{\E}_{1}, \dots, \wt{\E}_{n+1})$ forms an exceptional collection.
Denote it by $\wt{\sigma}.$

Let us check that the collection $\wt{\sigma}$ is strong and calculate the endomorphism
algebra of this collection.
Taking in account the exact sequence (\ref{projective}) we obtain the following isomorphisms
\[
\Hom^{j}_{\wX}(\pi^*\N,\; \wt{\E}_{n+1})\cong
\Hom^{j}_{\wX}(\pi^*\N,\; \pi^*\M)\cong
\Hom^{j}_{X}(\N,\; \M).
\]
If now $\N=u(N),$ where $N$ is an $\a$\!-module, then we obtain isomorphisms
\begin{equation}\label{homp}
\Hom^{j}_{\wX}(\pi^*\N,\; \wt{\E}_{n+1})\cong
\Hom^{j}_A(N,\; M).
\end{equation}
This implies that $\Ext^j_{\wX}(\wt{\E}_{i},\; \wt{\E}_{n+1})=0$ for all $i=1,\dots, n$ and $j\ge 1.$
Hence, the collection $\wt{\sigma}$ is strong. Moreover, the isomorphisms (\ref{homp}) for $j=0$
allow us calculate the endomorphism algebra of the collection $\wt{\sigma}.$ Thus,
the endomorphism algebra
$\End({\mathop{\bigoplus}\limits_{i=1}^{n+1}}\wt{\E}_i)$ of the collection $\wt{\sigma}$
is isomorphic to the ordinary extended algebra $\wa.$

Summarizing, we have an exact functor $\wt{u}: \mod-\wa\to \coh(\wX)$ that sends the
projective modules $\wt{P}_i$ to the vector bundles $\wt{\E}_i$ while simple modules $\wt{S}_i$ go to $\wt{\L}_i$ for all $i=1,\dots, n+1.$
Since $\End({\mathop{\bigoplus}\limits_{i=1}^{n+1}}\wt{\E}_i)\cong\wa$ and the collection $\wt{\sigma}$ is strong,
the derived
functor $\wt{\u}:\db{\mod-\wa}\to \db{\coh \wX}$ is fully faithful by Theorem \ref{keller}.
It is evident from
the definition of $\wt{\L}_i$ and $\wt{\E_i}$ for $i=1,\dots, n$ as pull backs of vector bundles
from $X$ that
the restriction of the functor $\wt{u}$ on $\mod-\a$ is isomorphic to $\pi^*\circ u.$
This implies that the conditions (1) and (2) of (V) hold for the functor $\wt{u}.$

Finally, we have to show that the condition (3) also holds for an appropriate line bundle $\wL'$ on
$\wX.$  Choosing $\wL'$ as a line bundle $\cO_{\wX}(1)\otimes \pi^*\R^{\otimes s},$ where
$\R$ is an ample line bundle on $X$ and $s$ is sufficiently large, we can guarantee that the condition (3) will hold.
Indeed, since the rank of $\F$ is greater than $2$
the line bundle $\wL'$ belongs to ${}^{\perp}\wt{\u}(\db{\mod-\wa}).$
Besides, we have isomorphisms
\[
H^j(\wX, \; \wt{\L}_i^{-1}\otimes \wt{\L}')\cong H^j(X, \L^{-1}_i\otimes\F^{\vee}\otimes\R^{\otimes s})\quad
\text{and}
\quad
H^j(\wX, \; \wt{\L}_{n+1}^{-1}\otimes \wt{\L}')\cong H^j(X, S^2(\F^{\vee})\otimes\R^{\otimes s})
\]
for the cohomology of $\wt{\L}_{n+1}^{-1}\otimes \wt{\L}'$ and of $\wt{\L}_i^{-1}\otimes \wt{\L}',$ when $i=1,\dots,n.$

Taking a sufficiently large $s$ we obtain vanishing of cohomology for $j>0$ and can guarantee that
all these bundles are generated by global sections  on $X$ and on $\wX.$
\end{dok}

Proposition \ref{step} as an induction step implies our Main Theorem.
\begin{theorem}\label{main}
Let $\a$ be a quiver algebra on $n$ ordered vertices. Then there exist a smooth projective variety $X$ and an exact functor
$u:\mod-\a\to\coh(X)$ such that the following conditions hold
\begin{enumerate}
\item[1)] the induced derived functor $\u:\db{\mod-\a}\to \db{\coh X}$ is fully faithful;
\item[2)] simple modules $S_i$ go to line bundles $\L_i$ on $X$ under $u;$
\item[3)] any $\a$\!-module $M$ goes to a vector bundle on $X;$
\item[4)] the variety $X$ is a tower of projective bundles and has a full exceptional collection.
\end{enumerate}
\end{theorem}
\begin{dok}
The proof proceeds by induction on $n.$ The base of induction is $n=1$ and $\a=\kk.$
In this case $X=\PP^1,$ the functor $u$ sends $\a$ to $\cO_{\PP^1},$ and $\L=\cO(1).$
The inductive step is Proposition \ref{step}. By construction, the variety $X$ is a tower of projective bundles and, hence,
 has a full exceptional collection.
\end{dok}
\begin{corollary}\label{exc_colllection_bundles}
Let $\dA$ be a small DG category over $\kk$ such that the homotopy category
$\T=\Ho(\dA)$ has a full exceptional collection
$
\T=\langle E_1,\dots, E_n\rangle.
$
Then there exist a smooth projective scheme $X$ and fully faithful functor $\r:\T\to \db{\coh X}$
such that the functor $\r$ sends the exceptional objects $E_i$ to vector bundles $\E_i$ on $X.$
\end{corollary}
It follows directly from Theorem \ref{main} taking in account Corollary \ref{exceptional_generator}.

\begin{corollary}
Let $\a$ be a quiver algebra on $n$ ordered vertices. Then there exist a smooth projective variety
$X$ and a vector bundle $\E$ on $X$ such that $\End_X(\E)=\a$ and $\Ext^j_X(\E, \E)=0$ for all $j\ne 0.$
Moreover, they can be chosen so that the rank of $\E$ is equal to the dimension of $A.$
\end{corollary}
\begin{dok}
It follows from Theorem \ref{main}. The vector bundles $\E$ is isomorphic to the direct sum
$\mathop{\bigoplus}\limits_{i=1}^n \E_i.$ Since simple modules go to line bundles the rank of $\E$ is equal to the dimension of $A$ over $\kk.$
\end{dok}

The quiver algebras on $n$ vertices is a particular case of a finite dimensional algebra.
\begin{conjecture}
For any finite dimensional algebra $\Lambda$ there exist a smooth projective variety $X$ and a vector bundle $\E$ on $X$ such that
$\End_X(\E)=\Lambda$ and $\Ext^j_X(\E, \E)=0$ for all $j\ge 1.$
\end{conjecture}

This conjecture looks reasonable because of the following statement proven in \cite{O_glue}.
\begin{theorem}\cite[Th.5.3]{O_glue}
Let $\Lambda$ be a finite dimensional algebra over $\kk.$ Assume that $S=\Lambda/\rd$ is a separable $\kk$\!--algebra.
Then there are a smooth projective scheme
$X$ and a perfect complex $\E^{\cdot}$ such that $\End(\E^{\cdot})\cong\Lambda$ and $\Hom(\E^{\cdot}, \E^{\cdot}[l])=0$
for all $l\ne 0.$
\end{theorem}

\section{Examples and applications}

As applications we consider two examples of realizations of quiver algebras as endomorphism algebras of
vector bundles on smooth projective varieties. 

\subsection{The quiver associated with the Ising 3-point function}
We will combine the methods of proof of Proposition \ref{step}
and Proposition \ref{simple} to obtain varieties of small dimensions.
Let us consider a quiver of type $(2,2;2)$ that is defined by the following rule
\begin{equation}\label{Ising}
Q_I=\Bigl(
\xymatrix{
\underset{1}{\bullet}\ar@<1ex>[rr]^{a_1}\ar@<-1ex>[rr]_{b_1} &&
\underset{2}{\bullet} \ar@<1ex>[rr]^{a_2} \ar@<-1ex>[rr]_{b_2} &&
\underset{3}{\bullet}
}
\Big|\quad  a_2 b_1=0,\; b_2 a_1=0
\Bigr).
\end{equation}

The compositions of arrows $a_2 a_1$ and $b_2 b_1$ give two arrows from the first vertex to the third that
will be denoted as $a$ and $b,$ respectively. In this case, the three projective modules $P_1, P_2, P_3$
form an exceptional collection
and all spaces of morphisms are 2-dimensional vector spaces.
We denote by $\wa$ the algebra of this quiver $Q_I.$ It is an ordinary extension of the algebra $\a$
of the quiver $\left( \bullet\rightrightarrows\bullet \right)$ generated by projective modules $P_1, P_2$ via the
 $\a$\!-module
$M=\Hom_{\wa}(P_1\oplus P_2, P_3).$

The quiver $Q_I$ corresponds to the Ising 3-point function and is related to
a Landau-Ginzburg model with superpotential the Weierstrass function $W(z) = \wp(z; \omega_1, \omega_2),$
and the identifications $z \sim z + n_1\omega_1 + n_2\omega_2,\; n_i \in \ZZ,$
where $\omega_i$ are the two periods of $\wp(z; \omega_1, \omega_2)$ (see e.g. \cite{CV}).
In other words, the quiver $Q_I$ in (\ref{Ising}) is a quiver of the category of D-branes of type A
in an LG model with the total space being an elliptic curve with a deleted flex point. Such a curve is given by an equation
$y^2=4x^2 - g_2x-g_3$ in the affine plane and is equipped with the superpotential $W(x,y)=x.$

Let us consider the two dimensional quadric $Y=\PP^1\times\PP^1$ and take the two line bundles $\cO_Y$ and $\cO_Y(2, -1)$ on it.
The pair $(\cO(2, -1), \cO)$ is exceptional and there is only nontrivial two-dimensional $\Ext^1$ from $\cO(2, -1)$
to $\cO.$ Let us consider the universal extension
\[
0\lto \cO^{\oplus 2}\lto \U\lto\cO(2, -1)\lto 0.
\]
As mutation in an exceptional pair the vector bundle $\U$ is exceptional too and the pair $(\cO, \U)$ is an exceptional pair.
Moreover, it is a strong exceptional pair, i.e. $H^j(Y, \U)=$ for $j>0$ and $H^0(Y, \U)=U$ is the two-dimensional vector space.

Let us fix two projective lines $L_1$ and $L_2$ on $Y$ that are fibers of the projection $Y=\PP^1\times \PP^1$ on the first component and consider the following short exact sequence
\[
0\lto \F\lto \cO(2, 0)^{\oplus 2}\stackrel{\phi}{\lto}\cO_{L_1}(1)\oplus\cO_{L_2}(1)\lto 0,
\]
where $\phi$ is a general morphism and $\F$ is the kernel of $\phi.$
The restriction of the line bundle $\cO(2, 0)$ on $L_i$ is the trivial line bundle. Since each sheaf $\cO_{L_i}(1)$ is generated by global sections
any general morphism $\phi$ is surjective and  consists of two components
$\phi_i: \cO(2, 0)^{\oplus 2}\to\cO_{L_i}(1)$ each of which is surjective.
It is easy to see that the map on global sections
\[
H^0(Y, \cO(2, 0)^{\oplus 2})\lto H^0(Y, \cO_{L_1}(1)\oplus\cO_{L_2}(1))
\]
is also surjective. Hence, $H^j(Y, \F)=0$ when $j>0$ and $H^0(Y, \F)$ is two-dimensional.
It can also be checked that $\Ext^j(\U, \F)=0$ for $j>0$ and
$\Hom(\U, \F)$ is two-dimensional for general morphism $\phi.$
Moreover, if we consider the functor $\bR\Hom^{\cdot}(\cO\oplus\U, -)$ from $\db{\coh Y}$ to $\db{\mod-\a},$
where $\a=\End_Y(\cO\oplus\U)=\End_{\wa}(P_1\oplus P_2),$ then this functor sends the bundle $\F$ to the $\a$-module
$M=\Hom_{\wa}(P_1\oplus P_2, P_3).$

Now we can apply Proposition \ref{simple} and consider the projective bundle $X=\PP(\F^{\vee})$ with the projection $\pi$ on $Y.$
By Proposition \ref{simple} the exceptional collection of vector bundles $\left( \cO_X, \pi^*\U, \cO_X(1)\right)$
is strong exceptional and the algebra of endomorphisms of this collection is the algebra $\wa$ that is the quiver algebra of the quiver $Q_I$
defined by the rule (\ref{Ising}).

The variety $X=\PP(\F^{\vee})$ is 3-dimensional smooth projective variety, that is rational and possesses a full exceptional collection.
The quiver $Q_I$ has an interesting property. It was checked many years ago by A.~Bondal that there is a module
 over the algebra of this quiver
that is exceptional but the semi-orthogonal complement to this module does not have any exceptional objects at all. Hence, this exceptional module
can not be include in a full exceptional collection in the derived category of modules over this algebra (see \cite{Kuz}).
This module as a representation of the quiver has $1$\!-dimensional vector spaces over each vertex and is such that
the
 $a$\!-arrows act as isomorphisms while the $b$\!-arrows act as $0.$ Thus, we obtain the following statement.

\begin{proposition}
There esist a smooth projective scheme $X$ that is a projectivization of two dimensional vector bundle over $\PP^1\times \PP^1$
and a strong exceptional collection of vector bundles $\left( \cO_X, \pi^*\U, \cO(1)\right)$ on it such that the algebra of endomorphisms of this collection
is exactly the algebra of the quiver $Q_I$ (\ref{Ising}). Moreover, the variety $X$ possesses a full exceptional collection (of length 8), on the one hand,
and, on the other hand, there is an exceptional collection of length 6 that can not be extended to a full exceptional collection on $X.$
\end{proposition}

\begin{remark}\label{remaark_line}{\rm
It is also useful to note that the algebra of the quiver $Q_I$ can not be realized as the endomorphism
algebra of a strong exceptional collection of line bundles on a variety. Indeed,  any morphism of line bundles on a smooth
irreducible projective scheme is an isomorphism at the generic point which contradicts the fact that $b_2 a_1=0.$
}
\end{remark}

\subsection{Quivers of noncommutative projective planes}

Noncommutative deformations of the projective plane can be described in terms
of exceptional collection \cite{ATV, BP}.
We know that the derived category $\db{\coh \PP^2}$ has a full strong exceptional collection of line bundles
 $(\cO, \cO(1), \cO(2)).$
This means that the category $\db{\coh \PP^2}$ is equivalent to the derived category of finite modules
of the quiver algebra for the following quiver with relations
\begin{equation}\label{pplane}
Q_{\PP^2}=\Bigl(
\xymatrix{
\underset{1}{\bullet}\ar@<2ex>[rr]^{a_1}\ar@<-2ex>[rr]_{c_1} \ar[rr]|{\;\; b_1\;} &&
\underset{2}{\bullet} \ar@<2ex>[rr]^{a_2} \ar@<-2ex>[rr]_{c_2} \ar[rr]|{\;\; b_2\;} &&
\underset{3}{\bullet}
}
\Big|\quad  a_2 b_1=b_2 a_1,\quad a_2 c_1= c_2 a_1,\quad b_2 c_1=c_2 b_1
\Bigr).
\end{equation}

A deformation of the category $\db{\coh \PP^2}$ is directly related to deformations of the relations of the quiver $Q_{\PP^2}.$
Namely, the derived category of coherent sheaves on a noncommutative projective plane
should be a triangulated category with a full strong exceptional collection
$(\F_0, \F_1, \F_2),$ for which
the spaces of homomorphisms from  $\F_i$ to $\F_j$ when $j-i=1$ are 3-dimensional
and the space $\Hom(\F_0, \F_2)$ is a 6-dimensional
vector space.
 Any such category  is determined by
a composition tensor $\mu: V\otimes U\to W,$ where $\dim V=\dim U=3$ and
$\dim W=6.$ This map should be surjective.
Thus, the derive category of coherent sheaves $\db{\coh \PP^2_{\mu}}$ on a noncommutative projective plane
$\PP^2_{\mu}$ is a category having a full strong exceptional collection
with composition tensor $\mu.$

Denote by $I$ the relations, i.e. the kernel of $\mu,$  and denote by $\nu$ the inclusion $I\to  V\otimes U.$
We will consider only the nondegenerate (geometric)
case, where the restrictions $\nu_{\bar{u}}: I\to V$ and
$\nu_{\bar{v}}: I\to U$ have rank at least two for all nonzero elements
$\bar{u}\in U^{\vee}$ and $\bar{v}\in V^{\vee}.$
The equations $\det \nu_{\bar{u}}=0$ and $\det\nu_{\bar{v}}=0$ define closed subschemes
$\Gamma_U\subset\PP(U^{\vee})$ and $\Gamma_V\subset\PP(V^{\vee}).$
It is easy to see that the correspondence which attaches  the kernel of the map
$\nu_{\bar{v}}^{\vee}: U^{\vee}\to I^{\vee}$ to
a vector $\bar{v}\in V^{\vee}$ defines
an isomorphism between $\Gamma_U$ and $\Gamma_V$.
Moreover, under these circumstances $\Gamma_U$ is either a cubic
in $\PP(U^{\vee})$ or the entire
projective plane $\PP(U^{\vee}).$ If $\Gamma_U=\PP(U^{\vee}),$
then $\mu$ is the standard tensor $U\otimes U\to S^2 U.$
In this case we obtain the usual projective plane $\PP^2.$

Thus, the non-trivial case is the situation where $\Gamma_V$ is a cubic. Let us denote it by $E$.
The curve $E$ comes equipped with two embeddings into the
projective planes $\PP(U^{\vee})$ and $\PP(V^{\vee}),$ respectively. The restriction of
the line bundles $\cO(1)$ these embeddings
determine two line bundles $\L_1$ and $\L_2$ of degree $3$ on $E.$
This construction has an inverse.

\begin{construction}\label{constr:mu}
{\rm
The tensor $\mu$ can be reconstructed from the triple
$(E, \L_1, \L_2).$
Namely, the spaces $U, V$ are isomorphic to
$H^0(E, \L_1)$ and $H^0(E, \L_2),$ respectively, and the tensor
$\mu: V\otimes U\to W$ is nothing but the canonical map
$
H^0(E, \L_2) \otimes H^0(E, \L_1)\lto H^0(E, \L_2\otimes\L_1).
$
}
\end{construction}

Note also that the mirror symmetry relation for
noncommutative planes is described in \cite{AKO} as a special elliptic fibration over
$\AA^1$ with three  ordinary critical points and with  symplectic forms, variations of which are
related to noncommutative deformations of $\PP^2.$

Let us fix a noncommutative projective plane $\PP_{\mu}$ that is defined by a
 tensor $\mu: V\otimes U\to W,$ where $\dim V=\dim U=3$ and
$\dim W=6.$
Consider the usual projective plane $\PP(U)$ and the vector bundle $\T(-1)$ on it.
The space of global section of $\T(-1)$ is canonically isomorphic to $U.$
The tenor $\mu$ defines the tensor $\nu: I\to  V\otimes U$ as above, where $I$ is the kernel of $\mu.$
They induce a morphism of vector bundles on $\PP(U)$
\[
\wt{\nu}: I\otimes \cO_{\PP^2}\lto V\otimes \T(-1),
\]
the cokernel of which is a 3-dimensional vector bundle on the projective plane $\PP(U)$ that we denote as $\F.$
Let us take the projective bundle $X=\PP(\F^{\vee})\stackrel{\pi}{\lto} \PP(U)$ and consider the tautological
line bundle $\cO_X(-1).$ Denote by $\L$ its dual
$\cO_X(1).$ There is an isomorphism $\bR \pi_*\L=\F.$
It was proved in Proposition \ref{simple} that the sequence $\sigma=(\cO_X, \pi^*\T(-1), \L)$
is strong exceptional. Moreover, it follows from the construction  that
\[
\begin{array}{l}
\Hom_X(\cO_X,\; \pi^*\T(-1))=\Hom_{\PP(U)}(\cO_{\PP^2},\; \T(-1))=U,\\
\Hom_X(\pi^*\T(-1),\; \L)=\Hom_{\PP(U)}(\T(-1),\; \F)=V,\\
\Hom_X(\cO_X,\; \L)=\Hom_{\PP(U)}(\cO_{\PP^2},\; \F)=W,
\end{array}
\]
and
the tensor of this collection is exactly $\mu: V\otimes U\to W.$
Hence, the subcategory of $\db{\coh X}$ generated by
the collection $\sigma$ is equivalent to the  derived category $\db{\coh \PP^2_{\mu}}$
of the noncommutative projective plane $\PP^2.$

Note that the derived category $\db{\coh X}$ has the following semi-orthogonal decomposition
\[
\db{\coh X}=\langle \pi^*\db{\coh \PP(U)}\otimes \L^{-1},\;
\pi^* \db{\coh \PP(U)},\; \pi^* \db{\coh \PP(U)}\otimes \L\rangle,
\]
where all three pieces are the derived categories of coherent sheaves of the usual projective plane.

On the other hand, since
$
\Hom_X(\pi^*\cO_{\PP^2}(1),\; \L)=\Hom_{\PP(U)}(\cO_{\PP^2}(1),\; \F)=0,
$
the line bundles $\pi^*\cO_{\PP^2}(1)$ and $\L$ are mutually orthogonal and we get the following decomposition
\[
\db{\coh X}=\left\langle \pi^*\db{\coh \PP(U)}\otimes \L^{-1},\;
\db{\coh \PP^2_{\mu}},\; \left\langle  \pi^*\cO_{\PP^2}(1),\; \pi^*\T(-1)\otimes\L,\; \pi^*\cO_{\PP^2}(1)\otimes \L \right\rangle\right\rangle.
\]
The subcollection $\sigma'=(\pi^*\cO_{\PP^2}(1),\; \pi^*\T(-1)\otimes\L,\; \pi^*\cO_{\PP^2}(1)\otimes \L)$ is also strong exceptional and there are
the following isomorphisms
\[
\begin{array}{l}
\Hom_Y(\pi^*\cO_{\PP^2}(1),\; \pi^*\T(-1)\otimes\L)=\Hom_{\PP(U)}(\T(-1),\; \F)=V, \\
\Hom_Y(\pi^*\T(-1)\otimes \L,\; \pi^*\cO_{\PP^2}(1)\otimes \L)=
\Hom_{\PP(U)}(\T(-1),\; \cO_{\PP^2}(1))=U,\\
\Hom_Y(\pi^*\cO_{\PP^2}(1),\; \pi^*\cO_{\PP^2}(1)\otimes \L)=\Hom_{\PP(U)}(\cO_{\PP^2},\; \F)=W.
\end{array}
\]

These isomorphisms show that the composition low in the collection $\sigma'$ is
given by the tensor $\mu^{\op}: U\otimes V\to W$ that is opposite to the tensor $\mu: V\otimes U\to W,$
i.e. $\mu^{\op}=\mu \cdot \iota,$ where $\iota: V\otimes U\to U\otimes V$ is the commutativity isomorphism.
Hence, the endomorphism algebra of the collection $\sigma'$ is opposite of the endomorphism algebra
of the collection $\sigma.$
We obtain the following semi-orthogonal decomposition for the derived category of coherent sheaves on $X$
\[
\db{\coh X}=\left\langle \db{\coh \PP(U)},\;
\db{\coh \PP^2_{\mu}},\; \db{\coh \PP^2_{\mu^{\op}}} \right\rangle.
\]

It should be noted that the noncommutative plane defined by the tensor $\mu^{\op}$
is related to the triple $(E, \L_2, \L_1).$ It can be easily checked that the triple
 $(E, \L_2, \L_1)$ is isomorphic to the original triple
$(E, \L_1, \L_2).$ Indeed, there is an automorphism $\tau$ of $E$ which is a multiplication by $-1$ with respect to
a suitable point of $E$ such that $\tau^* \L_1\cong \L_2$ and $\tau^* \L_2\cong \L_1.$
Thus, as an abstract noncommutative schemes $\PP^2_{\mu^{\op}}$ is isomorphic to
$\PP^2_{\mu}.$

Finally, let us mention one interesting fact about the variety $X.$ If
 $\mu$ is isomorphic to the usual tensor $U\otimes U\to S^2 U$ and we have the usual
 commutative projective plane $\PP(U^{\vee}),$ then the vector bundle $\F$ constructed above is isomorphic to
 the symmetric square
 $S^2(\T(-1)).$ In this case the variety $X=\PP(\F^{\vee})$ is isomorphic to the Hilbert scheme
 $\Hilb^2 \PP(U^{\vee})$
 of two points on the projective plane $\PP(U^{\vee}).$
 For a general tensor $\mu$ the variety $\PP(\F^{\vee})$ is a deformation of the Hilbert scheme
$\Hilb^2 \PP(U^{\vee}).$ Thus, any noncommutative plane can be obtained as an admissible subcategory
of a deformation of the Hilbert scheme of two points on the usual (dual) projective plane.
\begin{proposition}
For any noncommutative plane $\PP^2_{\mu}$ there is a 4-dimensional smooth projective variety $X$
of the form $\PP(\F^{\vee}),$ whose the derived category
$\db{\coh X}$ has a  semi-orthogonal decomposition
of the following form
\[
\db{\coh X}=\left\langle \db{\coh \PP^2},\;
\db{\coh \PP^2_{\mu}},\; \db{\coh \PP^2_{\mu^{\op}}} \right\rangle,
\]
where $\mu^{\op}$ is the opposite tensor to the tensor $\mu.$
Moreover, the variety $X$ is a deformation of the Hilbert scheme
$\Hilb^2 \PP^2$ of two points on the projective plane.
\end{proposition}

\end{document}